\documentclass[11pt]{article}
    
\usepackage{amsmath,amsfonts,amsthm,amscd,amssymb,graphicx}

\usepackage[notcite,notref]{showkeys}
\usepackage{hyperref}

\usepackage[utf8]{inputenc}

\usepackage{subfigure}
\usepackage{xcolor}

\numberwithin{equation}{section}

\usepackage{hyperref}

\def\bel{\begin{equation*} \begin{aligned}}
\def\eel{\end{aligned} \end{equation*}}
\def\beln{\begin{equation} \begin{aligned}}
\def\eeln{\end{aligned} \end{equation}}
\def\beq{\begin{equation}}
\def\eeq{\end{equation}}
\newcommand{\bea}{\begin{eqnarray}}
\newcommand{\eea}{\end{eqnarray}}

\newtheorem{theorem}{Theorem}[section]

\def\rit{{\mathbb{R}}}

\begin{document}

\title{Subcritical bifurcations of shear flows} 

\author{Dongfen Bian\footnotemark[1] \and Shouyi Dai\footnotemark[1]
\and Emmanuel Grenier\footnotemark[1]
}

\maketitle

\footnotetext[1]{School of Mathematics and Statistics, Beijing Institute of Technology, $100081$, Beijing, China. Emails:
biandongfen@bit.edu.cn and daishoui@outlook.com}

%\footnotetext[2]{Academy of Mathematics and System Science, Chinese Academy of Science, Beijing, China, Email:
%emmanuelgrenier@amss.ac.cn}

%\footnotetext[3]{Laboratoire J.A.Dieudonné, I.U.F., 
%Université Côte d’Azur, Parc Valrose, 06108 Nice Cedex 02, France}

%%%%%%%%%%%%%%%%%%%%%%%%

\subsubsection*{Abstract}

%%%%%%%%%%%%%%%%%%%%%%%%

It is well-known that shear flows in a strip or in the half plane are unstable for the incompressible Navier-Stokes
equations if the viscosity $\nu$ is small enough, provided the horizontal wave number $\alpha$
lies in a small interval, between the so called lower and upper marginal stability curves.
Moreover, a Hopf bifurcation occurs at the upper marginal stability curve.
In this article, for various shear flows, we give numerical evidences that this bifurcation is {\it subcritical}.

%%%%%%%%%%%%%%%%%%%%%%%%%%%%%%%%%%%%%%%%%%%%%%%%%%%%%%%%%%%%

\section{Introduction}

%%%%%%%%%%%%%%%%%%%%%%%%%%%%%%%%%%%%%%%%%%%%%%%%%%%%%%%%%%%%

In this paper, we address the classical question of the bifurcation of shear flows for the following
incompressible Navier-Stokes equations in the half space $\Omega = \rit \times \rit_+$
 \beq \label{NS1} 
\partial_t u^\nu + (u^\nu \cdot \nabla) u^\nu - \nu \Delta u^\nu + \nabla p^\nu = f^\nu,
\eeq
\beq \label{NS2}
\nabla \cdot u^\nu = 0 ,
\eeq
together with the Dirichlet boundary condition
\beq \label{NS3} 
u^\nu = 0 \quad \hbox{when} \quad y = 0.
\eeq
A shear flow is a stationary solution of (\ref{NS1},\ref{NS2},\ref{NS3}) of the form
$$
U(y) = (U_s(y),0), \qquad f^\nu = (- \nu \Delta U_s,0),
$$
where $U_s(y)$ is a smooth function, vanishing at $y = 0$ and converging at infinity to some constant $U_+$. 
In this article, we will study in details
the ``exponential flow"
$$
U_s(y) = 1 - e^{-y},
$$
We will also study the case of the strip $\Omega = \mathbb{R} \times (-1,+1)$.
In this case, $U_s(y)$ is a smooth function which vanishes at $\pm 1$,
and we will focus on flows of the form
$$
U_s(y) = 1 - y^{2 p}
$$
where $p = 1$, $2$ or $3$, the case $p = 1$ corresponding to the classical Poiseuille flow.

It is well known in physics 
that such shear flows are spectrally unstable (namely, the
corresponding linearized equations admit an exponentially growing solution) provided
the Reynolds number $Re = \nu^{-1}$ is large enough and provided
the horizontal wave number $\alpha$ of the perturbation lies in some interval $[\alpha_-(\nu),\alpha_+(\nu)]$
depending on the viscosity.
The curve $\alpha_-(\nu)$ (respectively, $\alpha_+(\nu)$)
is called the lower marginal stability curve (respectively, the upper marginal stability curve).
We in particular refer to \cite{Reid,Reid2,Schmidt} or to \cite{BG1,BDG,Guo} for detailed discussions
from the physical or mathematical points of view. 

More precisely, the linearized Navier-Stokes equations may be written under the form
\beq \label{LNS}
\partial_t u^\nu = \boldsymbol{L}_\nu u^\nu
\eeq
with
$$
\boldsymbol{L}_\nu u^\nu = \Pi \Bigl[ - (U \cdot \nabla) u^\nu - (u^\nu \cdot \nabla) U + \nu \Delta u^\nu \Bigr],
$$
where $\Pi$ is the Leray projection on divergence free vector fields.
We also define the bilinear symmetric operator $B$ by
$$
B(u,v) = - {1 \over 2} \Pi \Bigl[  (u \cdot \nabla) v +  (v \cdot \nabla) u \Bigr].
$$
We now take the Fourier transform in the $x$ variable, with dual Fourier variable $\alpha$,
which leads to
\beq \label{LNS1}
\partial_t u^\nu_\alpha = \boldsymbol{L}_{\nu,\alpha} u^\nu_\alpha
\eeq
where $u^\nu_\alpha$ is the Fourier transform of $u^\nu$ and $\boldsymbol{L}_{\nu,\alpha}$ the operator $\boldsymbol{L}_\nu$
restricted to the $\alpha$ Fourier component.

Following \cite{BG1,BDG,Reid,Guo,Reid2,Schmidt}, for $\nu$ small enough, there exists two smooth  
functions $\alpha_-(\nu)$ and $\alpha_+(\nu)$ such that $\boldsymbol{L}_{\nu,\alpha}$ has an eigenvalue $\lambda$ with
$\Re \lambda > 0$ if and only if $\alpha_-(\nu) < \alpha < \alpha_+(\nu)$.
Moreover, in this range, for a given $\alpha$, such an eigenvalue $\lambda$ is unique. It will be denoted by $\lambda(\alpha,\nu)$.
We will denote by 
$$
\zeta_{\alpha,\nu}(x,y) = \nabla^\perp [ e^{i \alpha x}  \psi_{\alpha,\nu}(y) ]
$$
a corresponding eigenvector.
Following \cite{BDG,Reid}, there exists constants $C_-$ and $C_+$ such that
\beq \label{asymptot1}
\alpha_-(\nu) \sim C_- \nu^{1/4}, 
\qquad
\alpha_+(\nu) \sim C_+ \nu^{1/6}
\eeq
as $\nu$ goes to $0$.
When $\alpha = \alpha_+(\nu)$, 
$$
\lambda(\alpha_+(\nu),\nu) = i \omega_0(\nu) = - i \alpha c(\nu)
$$ 
is purely imaginary.
The curves $\alpha_\pm(\nu)$ are defined for $\nu < 1/Re_c$
where $Re_c$ is the critical Reynolds number of the flow.

Let us focus on the upper marginal stability curve.
Let us fix $\nu < Re_c^{-1}$ and let us restrict ourselves to periodic functions, with period $2 \pi / \alpha$.
As $\alpha$ decreases and crosses $\alpha_+(\nu)$, two conjugated eigenvalues of $\boldsymbol{L}_\nu$
cross the real axis, these eigenvalues corresponding to the Fourier modes $\pm \alpha$.
In other words, as proved in \cite{Iooss}, $\boldsymbol{L}_\nu$ has a Hopf bifurcation at $\alpha_+(\nu)$. 
More precisely, in \cite{Iooss}, we proved the following Theorem.

\begin{theorem} \label{maintheo2} \cite{Iooss}
Let $\nu_0$ be small enough. Let $\alpha = \nu_+(\nu_0)$.
Let 
$$
\mu = \nu - \nu_0.
$$
Then, for $\mu$ in a right or left neighborhood of $0$, there exists a bifurcating time-periodic solution $u^\nu$ 
of  the Navier-Stokes equations which is a traveling wave function 
\beq \label{travel}
u^\nu(t,x,y) = V_{\varepsilon } \Bigl( x + {\omega \over \alpha} t ,y \Bigr) 
\eeq
where $V_\varepsilon$ is of the form
\begin{equation*}
V_{\varepsilon } = U(y) + \varepsilon V_1 + \varepsilon^2 V_2 +\mathcal{O}(\varepsilon ^{3}),
\end{equation*}
with
\begin{align*}
V_1(x,y) &= \zeta(x,y) + \overline{{\zeta(x,y) }}, \\
V_2(x,y) &= \hat V_{2,2}(y) e^{2i\alpha x }+\overline{\hat V_{2,2}(y)} e^{-2i\alpha x  } + V_{2,0},
\end{align*}
where $\zeta = \zeta_{\alpha,\nu}$, 
\begin{align*}
\mu (\varepsilon )  &= \varepsilon^{2}\mu_{2}+\mathcal{O}(\varepsilon ^{4}),
\\
\omega (\varepsilon ) &= \omega_0 %+\varepsilon ^{2}\omega _{2}
+\mathcal{O}(\varepsilon ^{2})
\end{align*}%
and where
\beq \label{inver}
e^{2 i \alpha x} \hat V_{2,2} = (2i\omega _{0}-\boldsymbol{L}_{\nu_0})^{-1}B( \zeta,\zeta), 
\eeq
\beq \label{inver2}
V_{2,0} =-2 \boldsymbol{L}_{0}^{-1}B(\zeta,\overline{\zeta}),
\eeq
\beq \label{expressionC}
C = \Bigl\langle 4 B \Bigl(\zeta, \boldsymbol{L}_{0}^{-1} B(\zeta,\overline{\zeta}) \Bigr)
- 2 B \Bigl(\bar \zeta, (2 i \omega_0 - \boldsymbol{L}_{\nu_0})^{-1} B(\zeta,\zeta) \Bigr), 
\zeta^\star \Bigr\rangle 
\eeq
$$
\mu _{2}
=  \langle \Delta \zeta,\zeta^{\ast} \rangle^{-1} ( C + i \omega_2 ).
$$
The bifurcation is supercritical if $\Re C > 0$ and subcritical if $\Re C < 0$. 
\end{theorem}

The sign of $\Re C$ is crucial to understand the dynamics of the flow near the upper marginal stability curve since, 
in the case of the strip, the amplitude $A(t)$ of any perturbation of the original shear flow can be described by its Landau equation
\beq \label{Landau}
{dA  \over dt} = \lambda A - C A |A|^2.
\eeq
In particular 
\beq \label{Landau2}
{d|A|^2 \over dt} = 2 \Re \lambda |A|^2 - 2\Re{C} |A|^4.
\eeq
Just below the upper marginal stability curve, $\Re \lambda > 0$.
Thus, 

\begin{itemize}
\item if $\Re C > 0$, the solution $u^\nu$ of Theorem \ref{maintheo2} exists for $\Re \lambda >0$ and is linearly and nonlinearly stable. 

\item if $\Re C < 0$, the solution $u^\nu$ exists for $\Re \lambda <0$ and is unstable.
In this case, we expect that small perturbations of the initial shear flow $U_s$ grow in the region where $\Re \lambda>0$, 
and reach a size of order $O(1)$, leading to turbulence.
\end{itemize}
Note that (\ref{Landau}) is not valid in the half space case where the essential spectrum of $\boldsymbol{L}_{\nu_0}$ touches $0$ \cite{Iooss}, however we give (\ref{Landau}) as a formal argument, which is be rigorous in the case of a flow in a bounded strip.

Unfortunately, up to now, there is no theoretical argument to assert the sign of $\Re C$ and we have to 
rely on numerical computations to determine whether the bifurcation is super-critical or sub-critical.
The aim of this article is precisely to compute the sign of  $\Re C$ near the upper marginal stability curve
and to prove that the bifurcation is {\it subcritical}
for various shear layer profiles.

\bigskip

\noindent{\it Numerical evidence:
for the following shear flows,
at the upper marginal stability curve, the Hopf bifurcation is sub-critical:

\begin{itemize}

\item the exponential shear flow $U_s(y) = 1 - e^{-y}$ in the half space,

%\item the Gaussian shear flow $U_s(y) = 1-e^{-y^2}$, {\color{red} to be done}

%\item the classical Blasius shear flow, {\color{red} to be done}

\item Poiseuille type flows in the strip $-1<y<1$, $U_s(y) = 1 - y^{2p}$ with $p = 1$ (genuine Poisseuille flow), $p = 2$ and $p = 3$.

\end{itemize}
}

Note that this result is already well-known in physics for the Poisseuille flow
(see for instance \cite{Schmidt,Sen_Venkateswarlu_1983}), but we have not found a similar result for the
exponential shear flow in the literature.

\medskip

Let us turn to the lower marginal stability curve. As we cross this curve, the real parts of
two complex conjugate eigenvalues of $\boldsymbol{L}_{\nu}$ become positive. We have an Hopf bifurcation
and Theorem \ref{maintheo2} remains valid provided $2 \alpha_-(\nu) > \alpha_+(\nu)$, namely provided
the second harmonic is stable.
We already know that $\alpha_-(\nu) \ll \alpha_+(\nu)$ as $\nu$ goes to $0$, thus there exists
a Reynolds number $Re_d$ such that
$$
2 \alpha_-(Re_d^{-1}) = \alpha_+(Re_d^{-1}).
$$
Direct numerical computations lead to the following statement.

\bigskip

\noindent{\it Numerical evidence:
for all the previous shear flows, there exists two Reynolds numbers $Re_s$ and $Re_d$ such that
$Re_c < Re_s < Re_d$ and such that

\begin{itemize}

\item if $Re_c \le Re < Re_d$, a Hopf bifucation occurs at $\nu = Re^{-1}$ and
$\alpha = \alpha_-(Re^{-1})$.
In particular, near the lower marginal stability curve,
the Navier-Stokes system admits a traveling wave solution of the form (\ref{travel}).

\item if $Re_c \le Re < Re_s$, this bifurcating traveling wave exists in the region $\Re \lambda <0$ and is unstable,

\item if $Re_s \le Re < Re_d$, this bifurcating traveling wave exists in the region $\Re \lambda >0$ and is stable.

\end{itemize}
}

\section{Principle of the computation}

%%%%%%%%%%%%%%%%%%%%%%%%%%%%%%%%%%%%%%%%%%%%%%%%%%%%%%%%%%%%

The formulas to compute $C$ are explicitly given in Theorem \ref{maintheo2}.
Let us detail them step by step in the half-space case. The starting point is a solution $\psi_{\alpha,\nu}(y)$ 
of the Orr-Sommerfeld equation 
\beq \label{OS}
(U_s - c)  (\partial_y^2 - \alpha^2) \psi_{\alpha,\nu} - U_s''  \psi_{\alpha,\nu}  
- { \nu \over i \alpha}   (\partial_y^2 - \alpha^2)^2 \psi_{\alpha,\nu} = 0
\eeq
with corresponding eigenvalue $\lambda(\alpha,\nu) = i \omega_0 = - i \alpha c$
and boundary conditions $\psi_{\alpha,\nu}(0) = \partial_y \psi_{\alpha,\nu}(0) =0$
and $\psi_{\alpha,\nu} \to 0$ at infinity.
Note that the adjoint Orr-Sommerfeld equation is
$$
(U_s - c) (\partial_y^2 - \alpha^2) \psi^\star_{\alpha,\nu}  + 2U_s' \partial_y \psi^\star_{\alpha,\nu} 
- \frac{i\nu}{\alpha}(\partial_y^2 - \alpha^2)^2 \psi^\star_{\alpha,\nu} = 0
$$
with the same boundary conditions.
As $\nu$ goes to $0$, in the half space case,
the lower and upper marginal stability curves satisfy \cite{BDG,Reid,Reid2}
 $$
\alpha_-(\nu) \sim 1.005 \nu^{1/4} {(U_s'(0))^{5/4} \over U_+^{3/2}} 
$$
and
$$
\alpha_+(\nu) \sim \Bigl({1 \over 2 \pi^2} {U_s'(0)^{11} \over U_1''(0)^2}\Bigr)^{1/6}\nu^{1/6},
$$
where $U_+$ is the limit of $U_s$ at infinity.
We do not reproduce here the corresponding equivalents of $c_\pm(\nu)$, the values of $c$ corresponding
to $\alpha_\pm(\nu)$ which can be easily found in \cite{BDG}.

The velocity $\zeta_{\alpha,\nu}$ is given by
$$
\zeta(x,y) =\zeta_{\alpha,\nu}(x,y) = \nabla^\perp [ e^{i \alpha x}  \psi_{\alpha,\nu}(y) ]
= e^{i \alpha x} \left( \begin{array}{c}
- \partial_y \psi_{\alpha,\nu} \\ 
 i \alpha \psi_{\alpha,\nu}
\end{array} \right)
$$
and its vorticity is
$$
\omega_{\alpha,\nu}(x,y) = \nabla \times \zeta_{\alpha,\nu} = e^{i \alpha x}
\Bigl(\partial_y^2 \psi_{\alpha,\nu} - \alpha^2 \psi_{\alpha,\nu} \Bigr) .
$$
Moreover
$$
\nabla \times B(\zeta,\zeta) 
= i \alpha e^{2 i \alpha x} \Bigl(\partial_y \psi_{\alpha,\nu} \partial_y^2 \psi_{\alpha,\nu} 
-  \psi_{\alpha,\nu} \partial_y^3 \psi_{\alpha,\nu} \Bigr) = e^{2 i\alpha x} B_0(y)
$$
and
$$
\nabla \times B(\zeta,\bar \zeta) 
= - \alpha \Im \Bigl( \partial_y \bar \psi_{\alpha,\nu} \partial_y^2  \psi_{\alpha,\nu}
+  \bar \psi_{\alpha,\nu} \partial_y^3 \psi_{\alpha,\nu} 
 \Bigr).
$$
We now turn to the resolution of (\ref{inver}). Let us introduce the stream function $\psi_{2,2}$ linked
to the velocity 
$$
\zeta_{2,2}(x,y) = e^{2i \alpha x} \hat V_{2,2} = \nabla^\perp [ e^{2 i \alpha x} \psi_{2,2}].
$$
Then, taking the curl of (\ref{inver}), we obtain the following Orr-Sommerfeld equation on $\psi_{2,2}$
$$
2i \alpha \Bigl[ (U_s - c) (\partial_y^2 - 4\alpha^2) \psi_{2,2} - U_s'' \psi_{2,2} \Bigr]
- \nu (\partial_y^2 - 4\alpha^2)^2 \psi_{2,2} = 
B_0.
$$
Similarly, we introduce $\psi_{2,0}$, the stream function of $V_{2,0}$, which leads to
$$
V_{2,0} = \nabla^{\perp} [ \psi_{2,0}(y)] = ( - \partial_y \psi_{2,0}, 0).
$$
The Orr-Sommerfeld equation (\ref{inver2}) on $\psi_{2,0}$ degenerates into
$$
- \nu \partial_y^4 \psi_{2,0} 
= 2\nabla \times B(\zeta,\bar \zeta) .
$$
We now compute $B(\zeta,V_{2,0})$ and $B(\bar \zeta, \zeta_{2,2})$
$$
B_1 = 
B(\zeta,V_{2,0}) = 
-{1 \over 2} \Pi e^{i\alpha x} 
\begin{pmatrix} 
i \alpha  \partial_y \psi(y) \partial_y \psi_0(y) - i \alpha \psi(y) \partial_y^2 \psi_0(y) \\ 
\alpha^2 \psi(y) \partial_y \psi_0(y)
\end{pmatrix}
$$
$$
B_2 = B(\bar \zeta, \zeta_{2,2}) = 
 -{1 \over 2} \Pi e^{i\alpha x} \begin{pmatrix}
i\alpha\big(  \partial_y \overline{\psi} \partial_y \psi_2 + \overline{\psi} \partial_y^2 \psi_2 - 2 \psi_2 \partial_y^2 \bar \psi \big) \\
3 \alpha^2 \big( 2 \partial_y \overline{\psi} \, \psi_2 + \overline{\psi} \partial_y \psi_2 \big)
\end{pmatrix}
$$
We then have
$$
C = - 2 \Bigl\langle B_1 + B_2, \zeta^\star \Bigr\rangle
$$
where the projection operator $\Pi$ disappears in the scalar product, and
$$
\zeta^\star = \nabla^\perp [ e^{i \alpha x}  \psi_{\alpha,\nu}^\star(y) ]
= e^{i \alpha x} \left( \begin{array}{c}
- \partial_y \psi_{\alpha,\nu}^\star \\ 
 i \alpha \psi_{\alpha,\nu}^\star
\end{array} \right).
$$
Note that we need to normalize $\zeta$ and $\zeta^\star$ such that $\langle \zeta, \zeta^\star \rangle = 1$.

%%%%%%%%%%%%%%%%%%%%%%%%%%%%%%%%%%%%%%%%%%%%%%%%%%

\section{Numerical results}

%%%%%%%%%%%%%%%%%%%%%%%%%%%%%%%%%%%%%%%%%%%%%%%%%

We follow the classical approach \cite{Schmidt} and use the Chebyshev polynomials
to find the eigenvalues and eigenvectors of the Orr-Sommerfeld equation and its adjoint, 
and to solve (\ref{inver}) and (\ref{inver2}). We validate our code by reproducing the numerical
values of the Table $1$ of \cite{Chen_Joseph_1973}.

\medskip

\begin{figure}[hbt]
 \centering{
 \includegraphics[width=6cm]{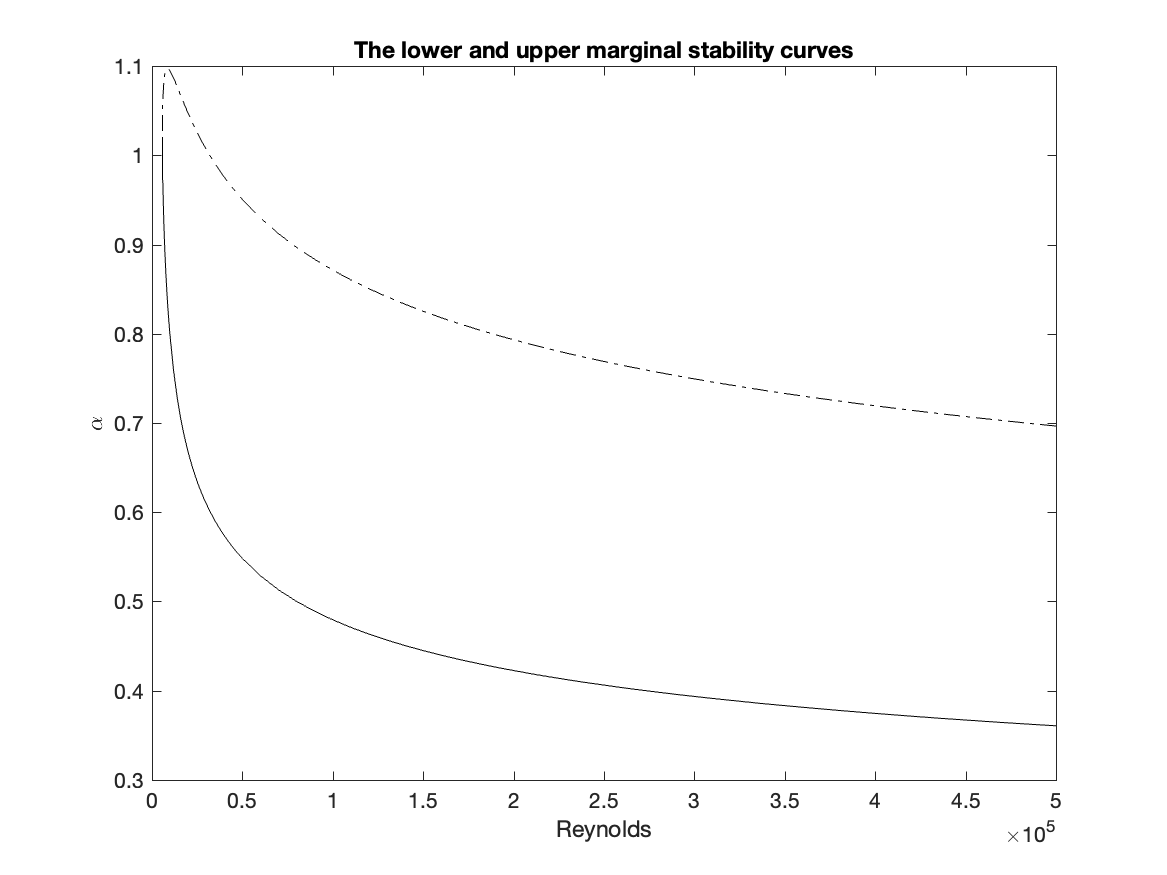}
 \includegraphics[width=6cm]{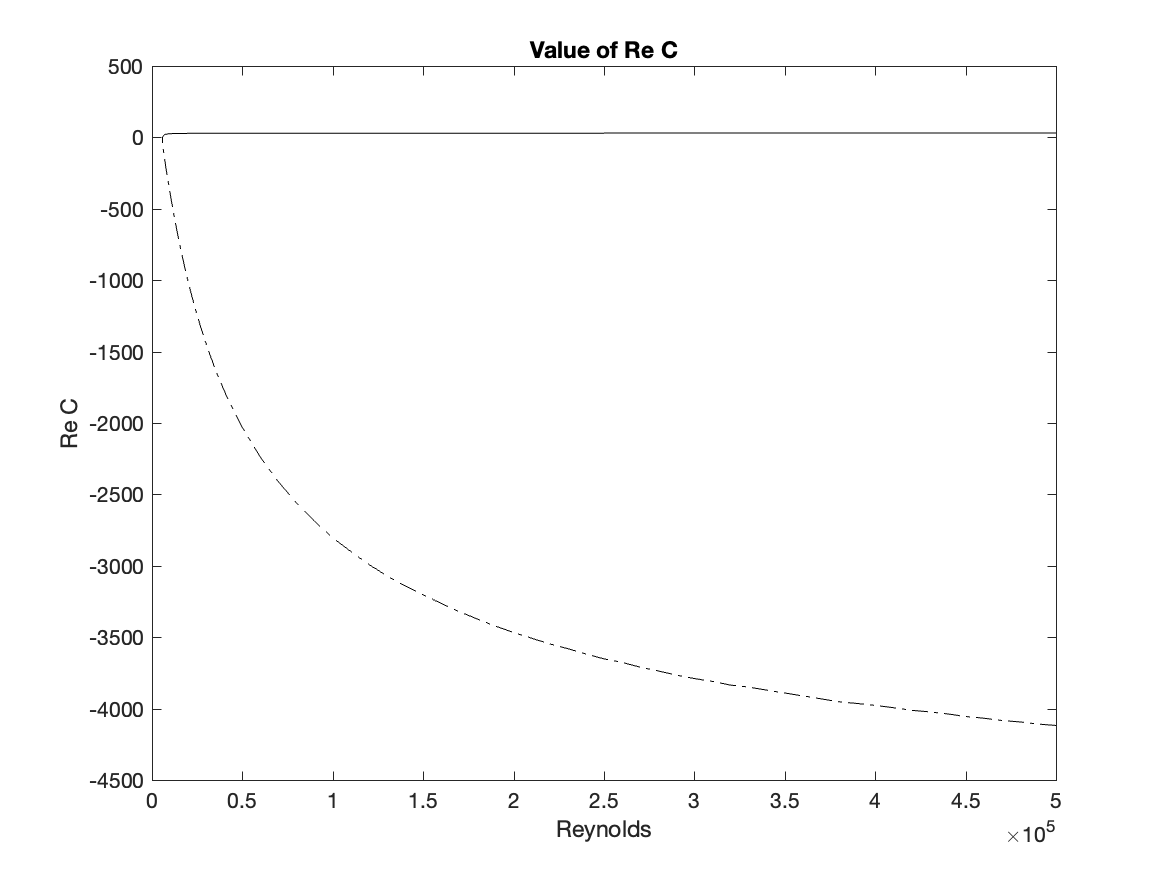}}
 \caption{
 Poiseuille}
 \label{figurePoiseuille}
\end{figure}

For Poiseuille flow, we recover Orzag's value of the critical Reynolds number $Re_c \approx 5772$.
The value of $\Re C$ is negative on the upper marginal stability curve.
It is negative on the lower one when $Re_c \le Re < Re_s \approx 5830$ and
 positive when $Re > Re_s$ (see Figure \ref{figurePoiseuille}).

\begin{figure}[hbt]
 \centering
 \includegraphics[width=6cm]{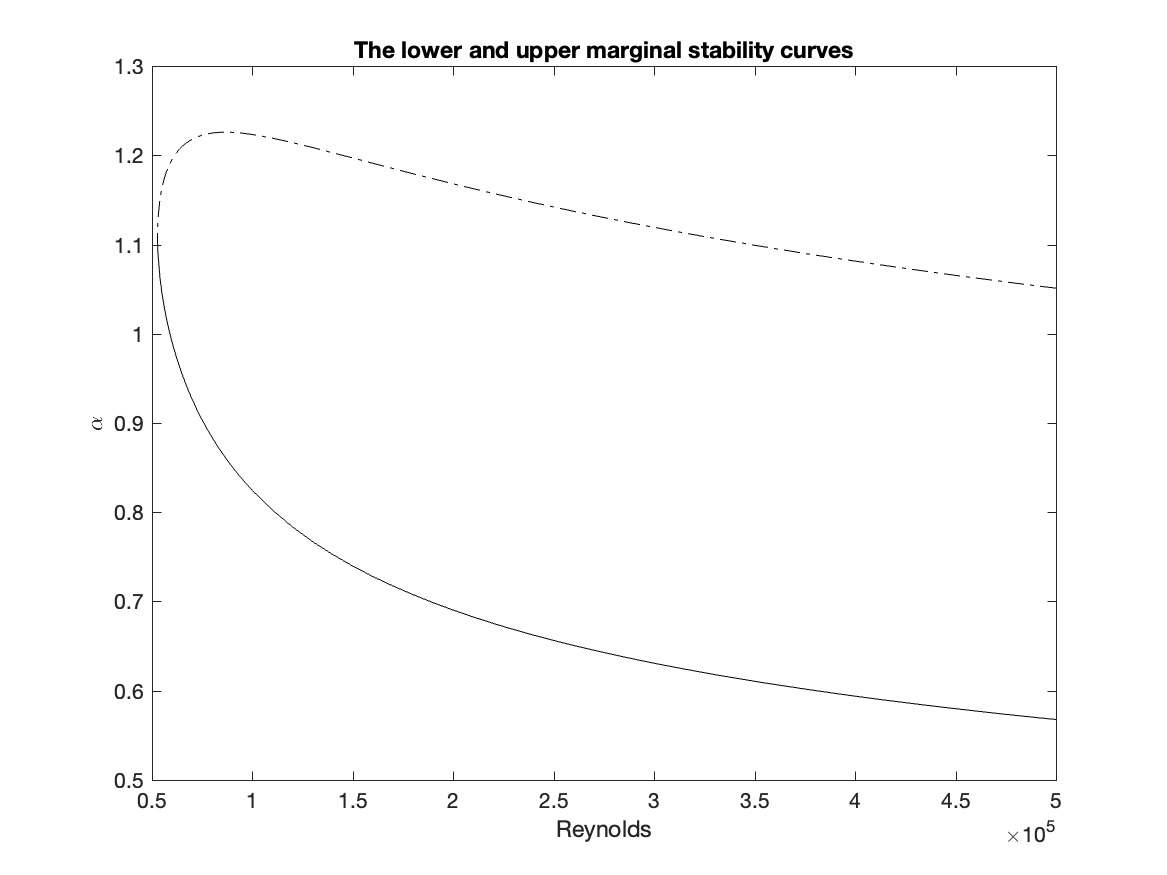}
 \hfill
 \includegraphics[width=6cm]{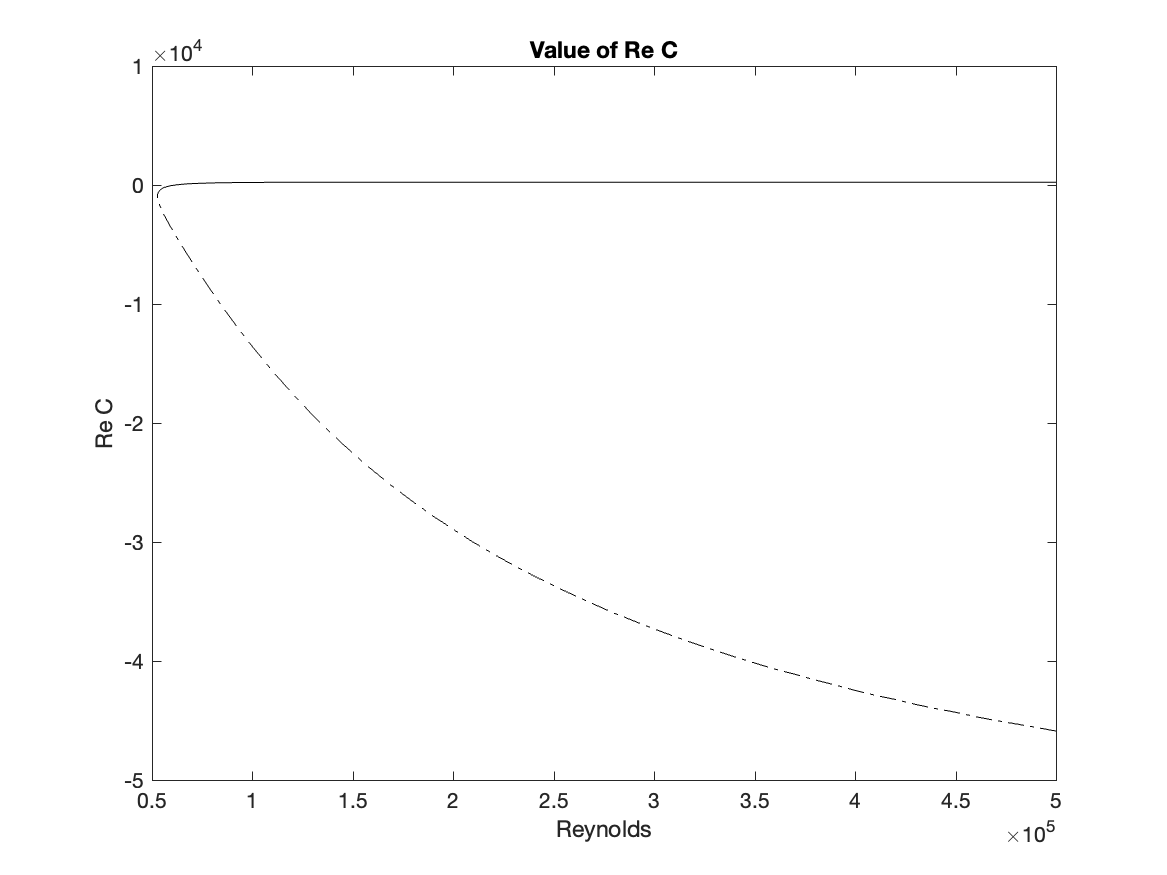}
 \caption{
 $U_s(y) = 1 - y^4$}
 \label{figureQuadric}
\end{figure}

\begin{figure}[hbt]
 \centering
 \includegraphics[width=6cm]{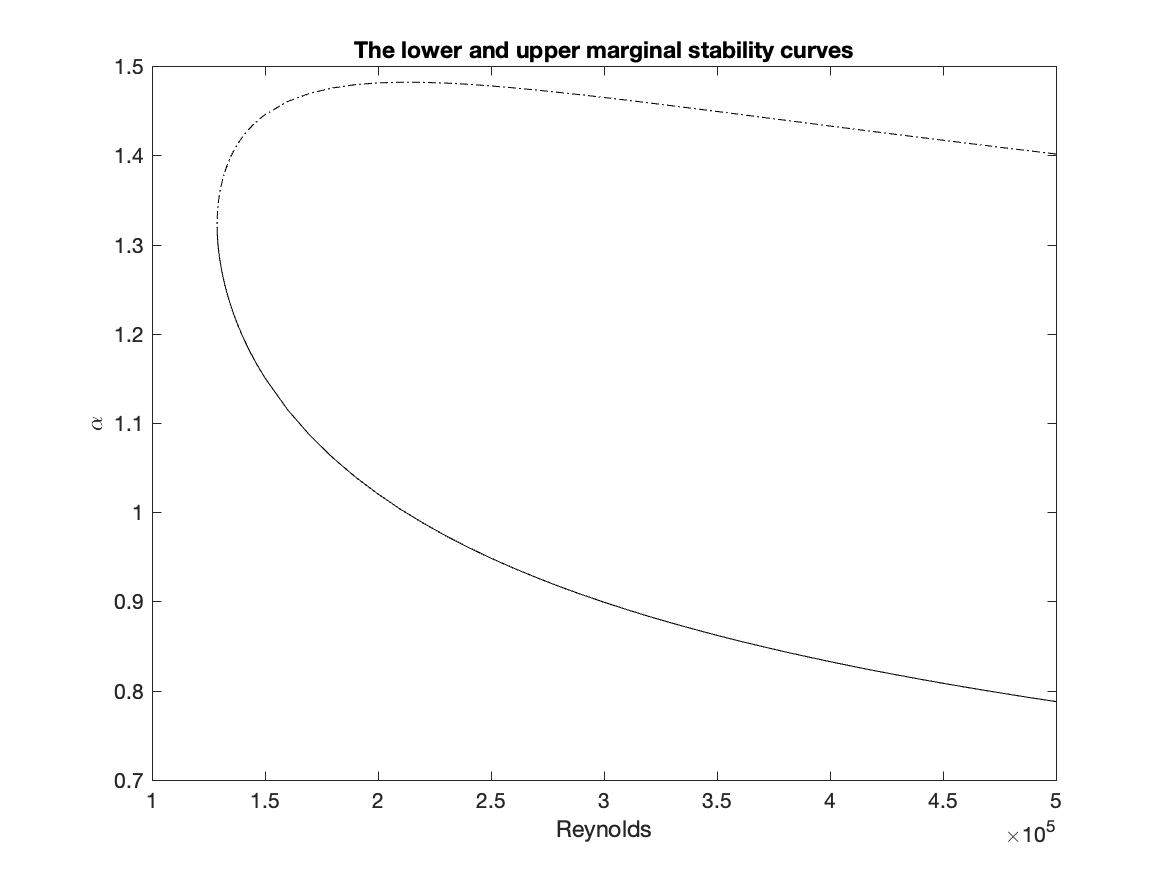}
 \hfill
 \includegraphics[width=6cm]{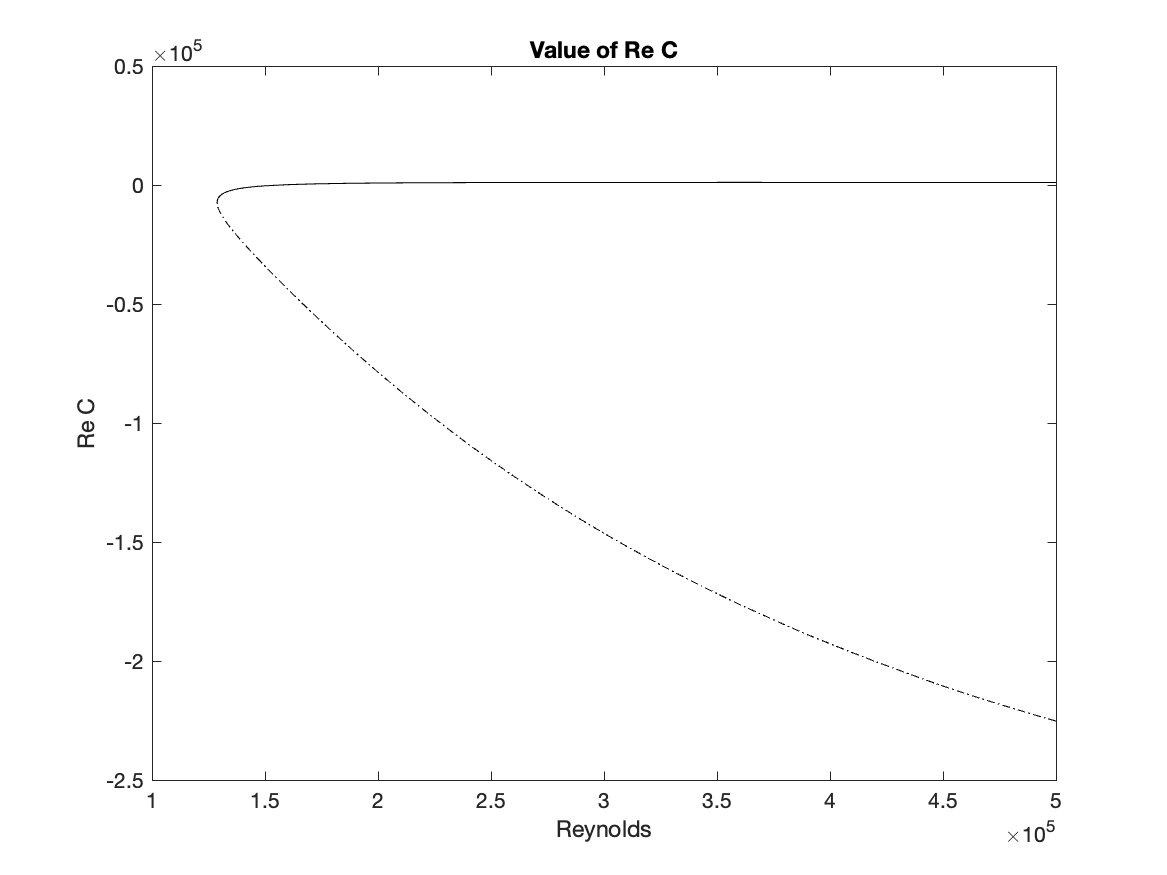}
 \caption{
 $U_s(y) = 1 - y^6$}
 \label{figureCub}
\end{figure}

The critical Reynolds number for $U_s(y) = 1 - y^{2p}$ for $p = 2$ (respectively $p = 3$)
is of order $52 748$ (respectively, $128820$). In these two cases, $\Re C < 0$ on the upper
marginal stability curve and positive on the lower one provided $Re > Re_s$ where $Re_s \approx 61461$ (respectively $156941$).
See Figures \ref{figureQuadric} and \ref{figureCub}.
For all the previous cases, $Re_d$ is very large, larger than $10^5$.

\medskip

Let us turn to the half space case. For the exponential profile $U_s(y) = 1 - e^{-y}$, we have $Re_c \approx 56375$.
The bifurcation is subcritical on the upper branch and also on the lower branch provided $Re < Re_s \approx 62714$.
Moreover, $Re_d \approx 85561$.

In all these cases, the bifurcation is {\it subcritical} on the upper marginal stability curve.

\begin{figure}[hbt]
 \centering
 \includegraphics[width=6cm]{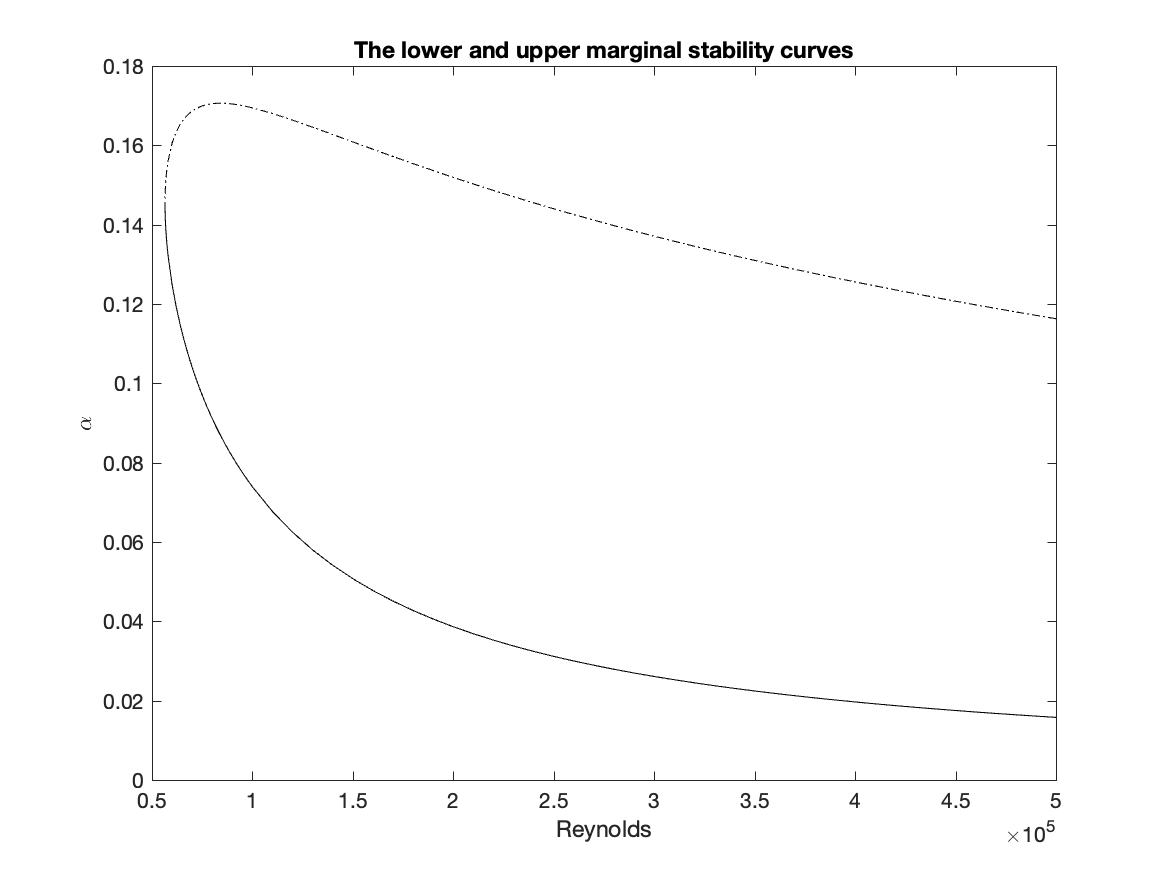}
 \hfill
 \includegraphics[width=6cm]{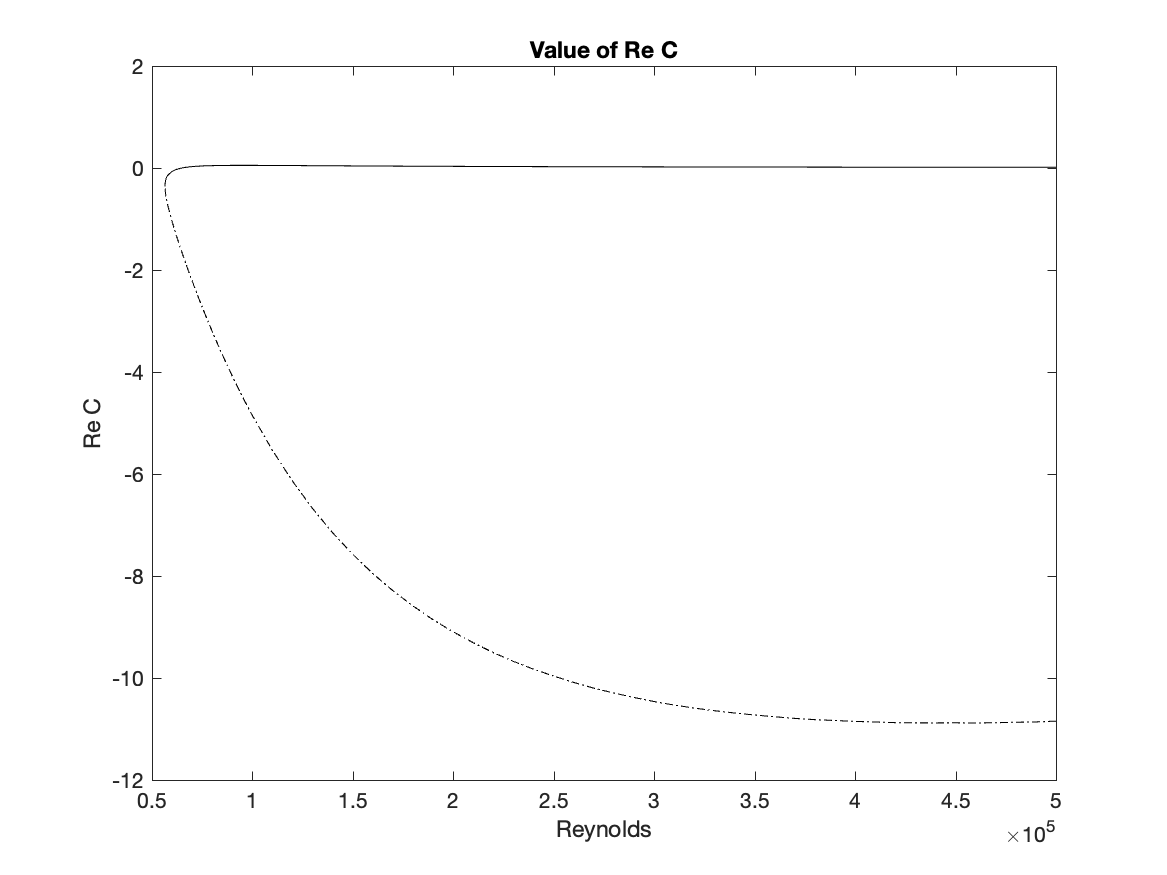}
 \caption{
 Exponential}
 \label{figureExponential}
\end{figure}

% We now turn to the second numerical approach, designed to handle very large Reynolds numbers.
% Figure \ref{figure1}  shows the dispersion relation, namely the imaginary part of $c$, and $Re C$ as a function of 
% $\alpha = \alpha_0 \nu^{1/6}$ when $Re = 10^{12}$ for the exponential profile.
% For this value of $Re$, the bifurcation is subcritical on the upper curve. 
% Moreover, $2 \alpha_-(\nu) < \alpha_+(\nu)$, thus Theorem \ref{maintheo2} is not valid in this area.

% \begin{figure}[hbt]
%  \centering
%  \includegraphics[width=1\textwidth]{Images2/dispersionRelation1d6.png}
%  \caption{
%  Dispersion Relation for $\alpha \sim \nu^{1/6}$  with $\nu =  10^{-12}$.
%  The vertical dashed lines indicates the zeros of $\Im c$, namely the lower and upper marginal stability curve.}
%  \label{figure1}
% \end{figure}

%  \begin{figure}[hbt]
%  \centering
%  \includegraphics[width=1\textwidth]{Images/nu1e-12kappa1d6DiffN1.png}
%  \caption{
% Values of $\Re C$ as a function of $\alpha$}\label{fig:Afig1DiffN}
% \end{figure}

%%%%%%%%%%%%%%%%%%%%%%%%%%%%%%%%%%%%%%%%%%%%%%%%%%%%%%%%%%%%%

\subsubsection*{Acknowledgements}

%%%%%%%%%%%%%%%%%%%%%%%%%%%%%%%%%%%%%%%%%%%%%%%%%%%%%%%%%%%%%

The authors warmly thank G. Iooss for his help on this work.
D. Bian is supported by NSFC under the contract 12271032.

%%%%%%%%%%%%%%%%%%%%%%%%%%%%%%%%%%%%%%%%%%%%%%%%%%%%%%%%%%%%%

\end{document}